\documentclass[a4paper,12pt]{article}
\usepackage{graphicx}
\usepackage{a4wide}
\usepackage[english]{babel}
\usepackage{amsfonts}
\usepackage{amsmath}
\usepackage{amssymb}
\usepackage{dsfont}
\newtheorem{lemma}{Lemma}

\newtheorem{theorem}{Theorem}

\allowdisplaybreaks

\date{}

\newcommand {\E} {\mathbb{E}}

\newcommand {\e} {\varepsilon}

\newcommand {\R} {\mathbb{R}}

\newcommand {\ve} {\varepsilon}

\numberwithin{equation}{section} \numberwithin{theorem}{section}
\numberwithin{lemma}{section} 
\numberwithin{corollary}{section}
\makeatletter\def\blfootnote{\xdef\@thefnmark{}\@footnotetext}\makeatother

\def\varep{\varepsilon}

\begin{document}
\title{\bf The Kadec-Pe\l czynski theorem in $L^p$, $1\le p<2$}
\author{I.\ Berkes\footnote{ Graz University of Technology, Institute of
Statistics, Kopernikusgasse 24, 8010 Graz, Austria.  \mbox{e-mail}: \texttt{berkes@tugraz.at}. Research
supported by FWF grant P24302-N18 and OTKA grant K 108615.} \, and R.\ Tichy\footnote{Graz University of Technology,
Institute of Mathematics A, Steyrergasse 30, 8010 Graz, Austria. \mbox{e-mail}: \texttt{tichy@tugraz.at}.
Research supported by FWF grant SFB F5510.}}

\maketitle \vskip0.5cm

\abstract{By a classical result of Kadec and Pe\l czynski (1962), every normalized weakly null sequence in $L^p$, $p>2$
contains a subsequence equivalent to the unit vector basis of $\ell^2$ or to the unit vector basis of $\ell^p$.
In this paper we investigate the case $1\le p<2$ and show that a necessary and sufficient condition for the first
alternative in the Kadec-Pelczynski theorem is that the limit random measure $\mu$ of the sequence satisfies
$\int_{\mathbb{R}} x^2 d\mu (x)\in L^{p/2}$.}

\section{Introduction}

Call two sequences $(x_n)$ and $(y_n)$ in a Banach space
$(B, \| \cdot \|)$ equivalent if there exists a constant $K>0$
such that
$$
K^{-1}\big\Vert \sum_{i=1}^n a_i x_i\big\Vert \le \big\Vert
\sum_{i=1}^n a_i y_i \big\Vert \le K\big\Vert \sum_{i=1}^n a_i
x_i\big\Vert
$$
for every $n\ge 1$ and every $(a_1,\ldots,a_n)\in \R^n$.
By a classical theorem of Kadec and
Pe{\l}czynski \cite{kape}, any normalized weakly null sequence
$(x_n)$ in $L^p(0, 1)$, $p>2$ has a subsequence equivalent to the
unit vector basis of $\ell^2$ or to the unit vector basis of
$\ell^p$. In the case when
$\{|x_n|^p, n\ge 1\}$ is uniformly integrable, the first alternative
holds, while if the functions $(x_n)$ have disjoint support, the second
alternative holds trivially. The general case follows via a subsequence
splitting argument as in \cite{kape}.

The purpose of the present paper is to investigate the case $1\le p<2$ and to
give a necessary and sufficient condition for the first alternative in the
Kadec-Pe\l czynski theorem. To formulate our result, we use probabilistic terminology.
Let $1\le p<2$ and let $(X_n)$ be a sequence of random variables defined on a probability
space $(\Omega, {\cal F}, P)$;
assume that $\{|X_n|^p, \, n\ge 1\}$ is uniformly integrable and $X_n\to 0$ weakly in $L^p$.
(This is meant as $\lim_{n\to\infty} \E (X_n Y)=0$ for all $Y\in L^q$
where $1/p+1/q=1$. To avoid  confusion with weak convergence of probability measures and distributions,
the latter will be called convergence in distribution and denoted by $\buildrel {\cal D}\over \longrightarrow$.)
Using the terminology of \cite{bero}, we call a sequence $(X_n)$ of random variables
{\it determining} if it has a limit distribution relative to any set $A$ in the
probability space with $P(A) > 0$, i.e.~for any $A\subset\Omega$
with $P(A) > 0$ there exists a distribution function $F_A$ such
that
$$\lim\limits_{n \to\infty} P(X_n \le t\mid A) = F_A(t)$$
for all continuity points $t$ of $F_A$. Here $P(\cdot |A)$ denotes conditional probability
given $A$. (This concept is the same as that of stable convergence, introduced in \cite{renyi}.) Since $\{|X_n|^p, \, n\ge 1\}$ is uniformly integrable, the sequence $(X_n)$ is tight and thus by
an extension of the Helly-Bray theorem (see e.g.\ \cite{bero}), it contains
a determining subsequence. Hence in the sequel we can assume, without loss of generality, that the sequence
$(X_n)$ itself is determining. As is shown in \cite{ald}, \cite{bero}, for any determining sequence $(X_n)$
there exists a random measure $\mu$ (i.e.\ a measurable map from  $(\Omega, {\cal F}, P)$ to
$(\cal M, \pi)$, where $\cal M$ is the set of probability measures on $\mathbb R$ and $\pi$ is the Prohorov
distance, see Section 3) such that for any $A$ with $P(A) > 0$ and any continuity point $t$ of $F_A$ we have
\begin{equation}\label{(4)}
F_A(t) =\E_A(\mu (-\infty, t]),
\end{equation}
where $\E_A$ denotes
conditional expectation given $A$. We call $\mu$ the {\it
limit random measure\/} of $(X_n)$.
We will prove the following result.

\begin{theorem}\label{th2}
Let $1\le p < 2$ and let $(X_n)$ be a determining sequence of
random variables such that
$\|X_n\|_p=1$ $(n=1, 2, \ldots)$, $\{|X_n|^p, \, n\ge 1\}$ is uniformly integrable and $X_n\to 0$
weakly in $L^p$. Let $\mu$ be the limit
random measure of $(X_n)$. Then there exists a subsequence $(X_{n_k})$ equivalent to the unit vector
basis of $\ell^2$ if and only if
\begin{equation}\label{maincond}
\int_{-\infty}^\infty x^2 d\mu(x) \in L^{p/2}.
\end{equation}
\end{theorem}

\bigskip
By assuming the uniform integrability of $|X_n|^p$, we exclude "spike" situations leading to a subsequence equivalent to
the unit vector basis of $\ell^p$ as in the Kadec-Pelczynski theorem. It is easily seen that
(\ref{maincond}) (and in fact $\int_{-\infty}^\infty x^2 d\mu (x)<\infty$ a.s.) imply
that for any $\delta>0$ there exists a set $A \subset \Omega$ with $P(A)\ge 1-\delta$ and a subsequence $(X_{n_k})$ such that
$$\sup_{k\ge 1} \int_A |X_{n_k}|^2 dP <\infty.$$
Thus the first alternative in the Kadec-Pe\l czynski theorem 'almost' implies bounded $L^2$ norms.

Call a sequence $(X_n)$ of random variables in $L^p$ {\it almost symmetric} if for any $\varepsilon>0$ there exists a $K=K(\varepsilon)$
such that for any $k\ge 1$, any indices $j_1>j_2>\ldots j_k \ge K$, any permutation $(\sigma (j_1), \ldots \sigma (j_k))$  of $(j_1, \ldots j_k)$
and any $(a_1, \ldots, a_k)\in \mathbb{R}^k$ we have
$$  (1-\e) \|\sum_{i=1}^k a_i X_{j_i}\|_p \le \|\sum_{i=1}^k a_i X_{\sigma (j_i)}\|_p     \le (1+\e) \|\sum_{i=1}^k a_i X_{j_i}\|_p\, .$$
Once in Theorem \ref{th2} we found a subsequence $(X_{n_k})$ equivalent to the unit vector basis of $\ell^2$, a result of Guerre \cite{gu}
implies the existence of a further  subsequence $(X_{m_k})$ of $(X_{n_k})$ which is
almost symmetric. Note that this conclusion also follows from the proof of
Theorem \ref{th2}. Guerre and Raynaud \cite{gura} also showed that for any $1\le p<q<2$ there exists a sequence $(X_n)$ in $L^p$, equivalent to the
unit vector basis of $\ell^q$, but not having an almost symmetric subsequence.
No characterization for the existence of almost symmetric subsequences of $(X_n)$ in terms of the limit random measure of $(X_n)$ or related quantities is known.

\bigskip
\section{Some lemmas}

\medskip
The necessity of the proof of Theorem \ref{th2} depends on a  general structure theorem for lacunary sequences proved in \cite{bepe} (see
Theorem 2 of \cite{bepe} and the definition preceding it); for the convenience of the reader we state it here as a lemma.

\begin{lemma} \label{lemma1}
Let $(X_n)$ be a determining sequence of r.v.'s and $(\ve_n)$ a
positive numerical sequence tending to 0. Then, if the underlying probability space is rich enough, there exists a
subsequence $(X_{m_k})$ and a sequence $(Y_k)$ of discrete r.v.'s
such that
\begin{equation}\label{1}
P\bigl( |X_{m_k} - Y_k| \geq \ve_k \bigr) \leq \ve_k \quad k =
1,2 \dots
\end{equation}
and for each $k > 1$ the atoms of the $\sigma$-field
$\sigma \{ Y_1, \dots, Y_{k - 1} \}$ can be divided into two
classes $\Gamma_1$ and $\Gamma_2$ such that

\noindent
{\rm (i)} \ $\sum_{B \in \Gamma_1} P(B) \leq \ve_k$;

\noindent
{\rm (ii)} \,  For any $B \in \Gamma_2$ there exist
$P_B$-independent r.v.'s $\{ Z^{(B)}_j, j = k, k + 1, \dots \}$
defined on $B$ with common distribution function $F_B$ such that
\begin{equation}
P_B \bigl( |Y_j - Z^{(B)}_j | \geq \ve_k \bigr) \leq \ve_k \quad
j = k, k + 1, \dots
\end{equation}
Here $F_B$ denotes the limit distribution of $(X_n)$ relative to
$B$ and $P_B$ denotes conditional probability given~$B$.
\end{lemma}

Note that, instead of (\ref{1}),  in Theorem 2 of \cite{bepe} the conclusion is
$\sum_{k=1}^\infty |X_{m_k}-Y_k|<\infty$ a.s., but after a further thinning,
(\ref{1}) will also hold. The phrase "the underlying probability space is rich enough"
is meant in Lemma \ref{lemma1} in the sense that on the underlying space there exists
a sequence of independent r.v.'s, uniformly distributed over $(0, 1)$ and also independent
of the sequence $(X_n)$. Clearly, this condition can be guaranteed by a suitable enlargement
of the probability space not changing the distribution of $(X_n)$ and $\mu$ and thus this
assumption can be assumed without loss of generality.

Lemma \ref{lemma1} means that every tight sequence of r.v.'s has a subsequence which can be closely
approximated by an exchangeable sequence having a very simple structure, namely which is i.i.d.\
on each set of a suitable partition of the probability space. This fact is an "effective" form
of the general subsequence principle of Aldous \cite{ald} (for a related result see Berkes and Rosenthal \cite{bero})
and reduces the studied
problem to the i.i.d.\ case which will be handled by the classical concentration technique of L\'evy \cite{le}, as improved by Esseen
\cite{es}.

\begin{lemma}\label{lemma2}
Let $X_1, X_2, \dots, X_n$ be i.i.d.\ random variables with
distribution function $F$ and put $S_{n}= X_1 + \cdots +
X_n$. Then for any $t > 0$ we have
\begin{equation}\label{2.5}
P\biggl(\biggl| S_n \biggr| \leq t \biggr) \leq A {t \over
\sqrt n} \biggl[ \int\limits_{|x| < t} x^2 dF(x) - 2 \biggl(
\int\limits_{|x| < t} x dF(x)\biggr)^2 \biggr]^{-1/2}
\end{equation}
provided the difference on the right-hand side is positive and
$\int\limits_{|x| < t} dF(x) \geq 1/2$. Here $A$ is an absolute constant.
\end{lemma}

\medskip
{\it Proof}. Let $F^*$ denote the distribution function obtained
from $F$ by symmetrization. The left hand side of (\ref{2.5}) is clearly bounded by $Q_{S_n} (2t)$, where $Q_{S_n} (\lambda)=\sup_x P(x\le S_n\le x+\lambda)$ is the concentration function of $S_n$. By a well known concentration function inequality of Esseen (see \cite{es}, formula (3.3)) we have
\begin{align}\label{esseen}
&Q_{S_n} (\lambda) \le An^{-1/2} \left(\lambda^{-2} \int_{|x|< \lambda} x^2 dF^*(x)+ \int_{|x|\ge \lambda} dF^*(x)\right)^{-1/2}\nonumber \\
&\le A\lambda \, n^{-1/2} \left(\int_{|x|< \lambda} x^2 dF^*(x)\right)^{-1/2}
\end{align}
for any $\lambda>0$, where $A$ is an absolute constant. Thus the left hand side of (\ref{2.5}) is bounded by the last expression in (\ref{esseen}) with $\lambda=2t$
and thus to prove (\ref{2.5}) it suffices to show that $\int\limits_{|x|< t} dF(x) \geq 1/2$ implies
\begin{equation}\label{2.6}
\int\limits_{|x| < 2t} x^2 dF^*(x) \geq \int\limits_{|x| < t} x^2 dF(x)
- 2 \biggl(\int\limits_{|x| < t} x dF(x) \biggr)^2.
\end{equation}
Let $\xi$ and $\eta$ be independent r.v.'s with distribution
function $F$, set
$$
C = \{ | \xi - \eta | < 2t \}, \quad D = \{ |\xi | < t,
|\eta| < t \}.
$$
Then
\begin{align*}
& \int\limits_{|x| < 2t} x^2 dF^*(x)
= \int\limits_C (\xi - \eta)^2 dP \geq \int\limits_D (\xi - \eta)^2 dP \\
&= 2 \int\limits_{|\xi| < t} \xi^2 dP \cdot P(|\eta | < t) - 2
\biggl( \int\limits_{|\xi| < t} \xi dP \biggr)^2 \geq \int\limits_{|\xi|< t} \xi^2 dP -
2 \biggl( \int\limits_{|\xi| < t} \xi dP\biggr)^2
\end{align*}
since $P(|\eta| < t) \geq 1/2$. Thus (\ref{2.6}) is valid.

\medskip

\begin{lemma} \label{lemma4}
Let $(X_n)$ be a determining sequence of r.v.'s with limit
random distribution function $F_\bullet$. Then for any set $A \subset \Omega$
with $P(A) > 0$ we have
\begin{equation}\label{2.7}
\E_A \biggl( \int\limits^{+\infty}_{-\infty} x^2 dF_\bullet (x) \biggr) =
\int\limits^{+\infty}_{-\infty} x^2 dF_A(x)
\end{equation}
in the sense that if one side is finite then the other side is
also finite and the two sides are equal. The statement remains
valid if in (\ref{2.7}) we replace the intervals of integration by
$(-t, t)$, provided $t$ and $-t$ are continuity points of $F_\Omega$.
\end{lemma}

\medskip
We used here the notation $F_\bullet$ to distinguish it from the ordinary limit
distribution function of $(X_n)$.

\medskip
{\it Proof.} Assume that $t$ and $-t$ are continuity points of $F_\Omega$. As observed in \cite[p.\ 482]{bero},
$t$ and $-t$ are continuity points of $F_\bullet$ with probability 1 (and hence also for $F_A$ for any $A\subset \Omega$ with $P(A)>0$)
and thus almost surely
$$
\int_{|x|< t} x^2 dF_\bullet (x) = -\left[ x^2 (1-F_\bullet (x)+F_\bullet (-x))\right]_{-t}^t +\int_{|x|< t} (1-F_\bullet (x)+F_\bullet (-x))2x dx $$
as it is seen by splitting the integral on the left hand side into subintegrals over $(-t,0)$ and $(0, t)$ (the integral over $\{0\}$ clearly equals 0) and using integration by parts. The same formula holds with $F_\bullet$ replaced by $F_A$. Integrating the last relation over $A\subset \Omega$ and using (\ref{(4)}) and Fubini's theorem, we get the validity of (\ref{2.7}) over $(-t, t)$. Letting $t\to\infty$ we get (\ref{2.7}) over $(-\infty, \infty)$.

For the following lemma (which is the key tool for the proof of the sufficiency part of
Theorem \ref{th2}) we need some definitions.
Given probability measures $\nu_n,\nu$ on the Borel sets of a separable
metric space $(S,d)$ we say that $\nu_n \buildrel {\cal D}\over
\longrightarrow \nu$ if
\begin{equation}\label{20}
\int_S f(x) \nu_n(dx) \longrightarrow \int_S f(x)
\nu (dx) \ \hbox{ as }\ n\to \infty
\end{equation}
for every bounded, real valued continuous function $f$ on $S$.
(For equivalent definitions and properties of this convergence see \cite{bill}).
(\ref{20}) is clearly equivalent to
\begin{equation}\label{21}
Ef(Z_n) \longrightarrow Ef(Z)
\end{equation}
where $Z_n,Z$ are r.v.'s valued in $(S,d)$ ({\it i.e.}, measurable maps
from some probability space to $(S,d)$) with distribution $\nu_n,\nu$.
A class $\cal G$ of real valued functions on $S$ is called {\it locally
equicontinuous} if for every $\varep >0$ and $ x\in S$ there is
a $\delta = \delta (\varep, x) >0$ such that $y\in S$, $d(x,
y) \le \delta$ imply $\vert f(x)-f(y)\vert \le \varep$ for
every $f\in {\cal G}$.

\begin{lemma} \label{raolemma}  (Ranga Rao \cite{rao})
Let $(S,d)$ be a separable
metric space and $\nu,\nu_n$ \ \ $(n=1,2,\ldots)$ probability measures
on the Borel sets of $(S,d)$ such that $\nu_n\buildrel {\cal D}\over
\longrightarrow \nu$. Let $\cal G$ be a class of real valued functions on
$(S,d)$ such that

\smallskip\noindent
(a) $\cal G$ is locally equicontinuous

\smallskip\noindent
(b) There exists a continuous function $g\ge 0$ on $S$ such that
$\vert f(x)\vert \le g(x)$ for all $f\in {\cal G}$ and $ x
\in S$ and
\begin{equation} \label{rao1}
\int_S g(x) \nu_n (dx) \longrightarrow \int_S g(x)\nu(dx)\ (<\infty)
\ \hbox{ as }\ n\to \infty.
\end{equation}
Then
\begin{equation} \label{rao2}
\int_S f(x)\nu_n(dx) \longrightarrow \int_S f(x)\nu(dx) \ \hbox{ as }
\ n\to \infty
\end{equation}
uniformly in $f\in {\cal G}$.\smallskip\par
\end{lemma}

\section{Proof of Theorem \ref{th2}}

Let $(\Omega ,{\cal F}, P)$ be the probability space of the $X_n$'s and
${\bf X} = (X_1,X_2,\ldots)$\ ; let further $\mu$ be the limit
random measure of $(X_n)$.  Let $(Y_n)$ be a sequence of r.v.'s on
$(\Omega ,{\cal F}, P)$ such that, given {\bf X} and $\mu$, the r.v.'s
$Y_1,Y_2,\ldots\ $ are conditionally i.i.d.\ with distribution $\mu$,
{\it i.e.},

\begin{equation}\label{18t}
P(Y_1\in A_1,\ldots,Y_k\in A_k\vert {\bf X},\mu) = \prod_
{i=1}^k P(Y_i\in A_i\vert {\bf X},\mu) \ \hbox{ a.s.}
\end{equation}
\begin{equation}\label{19t}
P(Y_j\in A\vert {\bf X},\mu )= \mu (A)\ \hbox{ a.s.}
\end{equation}
for any $j,k$ and Borel sets $A,A_1,\ldots,A_k$ on the real line.  Such
a sequence $(Y_n)$ always exists after a suitable enlargement of the
probability space (in fact $(Y_n)$ exists on $(\Omega ,{\cal F}, P)$ if
$(\Omega,{\cal F},P)$ is atomless over $\sigma ({\bf X},\mu)$, see the
vector-valued version of Theorem (1.5) of \cite{bero}; see also the remark
preceding Theorem (1.3) in \cite[p.\ 479]{bero}) or, alternatively,
the sequence $(X_n)$ can be redefined, without changing its distribution,
on a standard sequence space over which $(Y_n)$ can be defined, see
\cite[p.\ 72] {ald}.  Clearly, $(Y_n)$ is an exchangeable sequence; we call it
the {\it limit exchangeable sequence} of $(X_n)$. It is not hard to see (cf. \cite{ald},
\cite{bero}) that there exists a subsequence $(X_{n_k})$ such that for every $k\ge 1$ we have
\begin{equation}\label{6c}
(X_{n_{j_1}},\ldots,X_{n_{j_k}}) \buildrel {\cal D}\over \longrightarrow
(Y_1,\ldots ,Y_k) \ \hbox{ if }\ j_1<\cdots <j_k\ ,\ j_1\to \infty.
\end{equation}
Note that the existence of a subsequence $(X_{n_k})$ and exchangeable $(Y_k)$ satisfying (\ref{6c}) was first proved by
Dacunha-Castelle and Krivine \cite{dck} via ultrafilter techniques.
The limit exchangeable sequence, as defined above, also has the following simple property, proved in \cite[Lemma 12]{ald}.

\begin{lemma} \label{aldlemma}
For every $\sigma ({\bf X})$-measurable {\rm r.v.}\
$Z$ and any $j\ge 1$ we have
$$(X_n,Z) \buildrel {\cal D}\over \longrightarrow (Y_j,Z)$$
\end{lemma}

As before, let $\cal M$ denote the set of all probability measures
on $\mathbb R$ and let $\pi$ be the Prohorov metric on $\cal M$
defined by
\begin{align*}\pi (\nu ,\lambda ) &= \inf \bigl\{ \varep > 0 : \nu (A)
\le \lambda (A^\varep )+\varep \ \hbox{ and} \cr
&\qquad \lambda (A) \le \nu (A^\varep) +\varep \ \hbox{ for all Borel
sets } \ A\subset \mathbb{R} \bigr\}.
\end{align*}
Here
$$A^\varep=\{x \in \mathbb{R}: |x-y|<\varep \ \text{for some} \ y\in A\}$$
denotes the open $\varep$-neighborhood of $A$. Let
\begin{equation}\label{e16}
S = \left\{ \nu \in {\cal M} :\int x d\nu (x) = 0\ ,\ \int x^2
d\nu (x) <+\infty \right\}.
\end{equation}
Since $\int_{-\infty}^\infty x^2 d\mu(x)<\infty$ a.s. (which
follows from (\ref{maincond})) and $\int_{-\infty}^\infty x
d\mu(x)=0$ a.s.\ by $X_n\to 0$ weakly,
we have
\begin{equation}\label{e7}
P\bigl\{\mu \in S\bigr\} = 1.
\end{equation}
Following Aldous \cite{ald} we define another metric $d$ on $S$ by
\begin{equation}\label{e17}
d(\nu ,\lambda) = \left( \int_0^1 \bigl( F_\nu^{-1} (x) -
F_\lambda^{-1} (x)\bigr)^2 \, dx\right)^{1/2}
\end{equation}
where $F_\nu$ and $F_\lambda$ are the distribution functions of
$\nu$ and $\lambda$, respectively, and $F^{-1}$ is defined by
$$F^{-1}(x) = \inf \bigl\{ t: \ F(t) \ge x\bigr\}, \qquad 0<x<1$$
for any distribution function $F$.  The right side of (\ref{e17})
equals $\Vert F_\nu^{-1} (\eta) - F_\lambda^{-1} (\eta )\Vert_2$
where $\eta$ is a random variable uniformly distributed in
$(0,1)$.  Since $F_\nu^{-1} (\eta)$ and $F_\lambda^{-1}(\eta)$ are
r.v.'s with distribution $\nu$ and $\lambda$, respectively (and
thus square integrable), it follows that $d$ is a metric on $S$.
It is easily seen (cf.\ \cite[p.\ 80]{ald} and relation (5.15) in
\cite[p.\ 74]{ald}) that $d$ is separable and  generates the same Borel
$\sigma$-field as $\pi$.
By the definition of $d$ we have, letting $0$ denote the zero
distribution,
\begin{equation} \label{ed}
\E d(\mu,0)^p = \E ( \text{Var} (\mu))^{p/2}= \E
\left(\int_{-\infty}^\infty  x^2 d\mu (x)\right)^{p/2}<\infty
\end{equation}
by our assumption (\ref{maincond}).
The following lemma expresses the crucial equicontinuity property of $d$.

\begin{lemma} \label{equicont}
Let \begin{equation}\label{e11} \psi (a_1,\ldots,a_n) = \big\Vert
\sum_{i=1}^n a_iY_i\big\Vert_p.
\end{equation}
Then we have
\begin{equation}\label{e10}
\left| {\Vert t+ \sum\limits_{k=1}^n a_k\xi_k^{(\nu)}\Vert_p -
\Vert t+\sum\limits_{k=1}^n a_k\xi_k^{(\lambda)} \Vert_p}
 \right| \le K d (\nu, \lambda) \psi (a_1,\ldots,a_n)
\end{equation}
for some constant $K>0$, every $n\ge 1$, $\nu,\lambda \in S$, real
numbers $t,a_1,\ldots,a_n$ and i.i.d.\ sequences
$(\xi_n^{(\nu)}),(\xi_n^{(\lambda)})$ with respective
distributions $\nu$ and $\lambda$.
\end{lemma}

Relation (\ref{e10}) means that the class of functions $\{
f_{t,{\bf a},n}\}$ defined by
\begin{equation}\label{e12}
f_{t,{\bf a},n} (\nu) = \psi ({\bf a})^{-1} \big\Vert
t+\sum_{k=1}^n a_k \xi_k^{(\nu)} \big\Vert_p \qquad {\bf a} =
(a_1,\ldots a_n) \ne {\bf 0}
\end{equation}
(where the variable is $\nu$ and $t,{\bf a},n$ are parameters) is
equicontinuous. In the context of unconditional convergence of
lacunary series, the importance of such equicontinuity conditions
was discovered by Aldous \cite{ald}. A similar condition in terms
of the compactness of the 1-conic class belonging to the type determined
by $(X_n)$ was given by Krivine and Maurey (see Proposition 3 in Guerre
\cite{gu}). The proof of our results is, however, purely
probabilistic and we will not use types.

\medskip
{\it Proof of Lemma \ref{equicont}.} We start with recalling the well known fact that if
$(\xi_n)$ is an i.i.d.\ sequence with $E\xi_n=0$, $E\xi_n^2
<+\infty$ then
\begin{equation}\label{e15a}
C\Vert \xi\Vert_1\left( \sum_{i=1}^k a_i^2\right)^{1/2} \le
\big\Vert \sum_{i=1}^k a_i\xi_i\big\Vert_p \le \Vert \xi\Vert_2
\left( \sum_{i=1}^k a_i^2\right)^{1/2}\quad 
\end{equation}
for any $1\le p<2$ and any $(a_1, \ldots, a_n)\in {\mathbb R}^n$,
where $C>0$ is an absolute constant. Since the $L^p$ norm of
$\sum_{i=1}^k a_i \xi_i$ in  (\ref{e15a}) cannot exceed the $L^2$
norm, the upper bound in (\ref{e15a}) is obvious, while the lower
bound is classical, see \cite{mz}.  Since $ \E\left|\sum_{i=1}^n
a_i Y_i \right|^p $ can be obtained by integrating $
\E\left|\sum_{i=1}^n a_i \xi_i^{(\omega)} \right|^p $ over
$\Omega$ with respect to $dP(\omega)$ where for each $\omega\in
\Omega$ the $\xi_i^{(\omega)}$ are i.i.d.\ with distribution
$\mu(\omega)$, relation (\ref{e15a}) implies that
\begin{equation}\label{ynorm}
A\left( \sum_{i=1}^k a_i^2\right)^{1/2} \le \big\Vert \sum_{i=1}^k
a_iY_i\big\Vert_p \le B\left( \sum_{i=1}^k a_i^2\right)^{1/2}
\end{equation}
where
$$ A= C\left[\E \left( \int_{-\infty}^\infty |x| d\mu (x)\right)^p\right]^{1/p}, \qquad B=
\left[\E \left( \int_{-\infty}^\infty x^2 d\mu
(x)\right)^{p/2}\right]^{1/p}. $$ By (\ref{maincond}) and since
the assumptions of Theorem \ref{th2} imply that $\mu$ is not concentrated at zero a.s., we have $0<A\le B<\infty$.

Turning to the proof of (\ref{e10}), note that the $L_p$ norms on
the left hand side depend of the sequences $(\xi_n^{(\nu)}),
(\xi_n^{(\lambda)})$ only through their distributions $\nu,
\lambda$, but not the actual choice of these i.i.d.\ sequences and
thus it suffices to verify (\ref{e10}) for any specific
construction. Let $(\eta_n)$ be a sequence of independent r.v.'s,
uniformly distributed over $(0,1)$.  Then $\xi_n^{(\nu)} =
F_\nu^{-1}(\eta_n)$ and $\xi_n^{(\lambda)} =
F_\lambda^{-1}(\eta_n)$ are i.i.d.\ sequences with distribution
$\nu$ and $\lambda$, respectively.  Using these sequences in
(\ref{e10}), the left hand side is at most $\Vert \sum_{i=1}^n a_i
(\xi_i^{(\nu)} - \xi_i^{(\lambda)})\Vert_p$ and since
$\xi_n^{(\nu)} - \xi_n^{(\lambda)} = F_\nu^{-1} (\eta_n) -
F_\lambda^{-1}(\eta_n)$ is also an i.i.d.\ sequence with mean 0
and variance $d(\nu, \lambda)^2$, using (\ref{e15a}) and the first
relation of (\ref{ynorm}) we get that the left hand side of
(\ref{e10}) is at most $K d(\nu,\lambda) \psi (a_1, \ldots, a_n)$
with some constant $K>0$. This completes the proof of Lemma \ref{equicont}.

With the equicontinuity statement of Lemma \ref{equicont} at hand, we can prove the sufficiency part of Theorem \ref{th2}
with a selection procedure similar to \cite{be89}. Assume that
$(X_n)$ satisfies (\ref{maincond}) and fix $0<\varep \le 1/2$.  We
shall construct a sequence $n_1<n_2<\cdots \ $ of integers such
that
\begin{equation}\label{24t}
(1-\varep) \psi (a_1,\ldots,a_k) \le \big\Vert \sum_{i=1}^k a_i X_{n_i}
\big\Vert_p \le (1+\varep)\psi (a_1,\ldots,a_k)
\end{equation}
for every $k\ge 1$ and $(a_1,\ldots,a_k)\in R^k$.
In view of (\ref{ynorm}), this will imply that $(X_{n_k})$ is equivalent to the unit vector base of $\ell^2$, but
it actually shows more, namely that under the assumptions of Theorem \ref{th2} there is a subsequence $(1+\varepsilon)$-equivalent
to the limit exchangeable sequence and hence $(1+\varepsilon)$-symmetric.

To construct $n_1$ we set
\begin{align*}
Q({\bf a},n,\ell) &= \vert a_1X_n + a_2Y_2+\cdots + a_\ell
Y_\ell \vert^p\cr
R({\bf a},\ell) &= \vert a_1Y_1+a_2Y_2+\cdots + a_\ell Y_\ell \vert^p\cr
\end{align*}
for every $n\ge 1$, $\ell \ge 2$ and ${\bf a} = (a_1,\ldots,a_\ell) \in
\mathbb{R}^\ell$.  We show that
\begin{equation}\label{25t}
E\left\{ {{Q({\bf a},n,\ell)}\over {\psi ({\bf a})^p}}
\right\} \longrightarrow  E\left\{ {{R({\bf a},\ell)}\over {\psi ({\bf a})^p}}
\right\} \ \hbox{ as }\ n\to \infty \quad \hbox{uniformly in }\ {\bf a},\ell.
\end{equation}
(The right side of (\ref{25t}) equals 1.)  To this end we recall that, given
{\bf X} and $\mu$, the r.v.'s $Y_1,Y_2,\ldots \ $ are conditionally
i.i.d.\ with common conditional distribution $\mu$ and thus, given
${\bf X},\mu$ and $Y_1$, the r.v.'s $Y_2,Y_3,\ldots \ $ are conditionally
i.i.d.\ with distribution $\mu$.  Thus
\begin{equation}\label{26t}
E\bigl( Q({\bf a},n,\ell)\vert {\bf X},\mu\bigr) =g^{{\bf a},\ell}
(X_n,\mu)
\end{equation}
and
\begin{equation}\label{27t}
E\bigl( R({\bf a},\ell)\vert {\bf X},\mu,Y_1\bigr) = g^{{\bf a},\ell}
(Y_1,\mu)
\end{equation}
where
$$g^{{\bf a},\ell} (t,\nu) = E\big\vert a_1t+\sum_{i=2}^\ell a_i\xi_i^
{(\nu)} \big\vert^p \qquad (t\in \mathbb{R}^1\ ,\ \nu \in S)$$
and $(\xi_n^{(\nu)})$ is an i.i.d.\ sequence with distribution $\nu$.
Integrating (\ref{26t}) and (\ref{27t}) we get
\begin{equation}\label{26at}
E\bigl( Q({\bf a}, n,\ell)\bigr) = Eg^{{\bf a},\ell}
(X_n,\mu)
\end{equation}
\begin{equation}\label{27at}
E\bigl( R({\bf a},\ell)\bigr) = Eg^{{\bf a},\ell} (Y_1,\mu)
\end{equation}
and thus (\ref{25t}) is equivalent to
\begin{equation}\label{28t}
E {{ g^{{\bf a},\ell} (X_n,\mu)}\over {\psi ({\bf a})^p}}
\longrightarrow
E {{g^{{\bf a},\ell} (Y_1,\mu)} \over {\psi ({\bf a})^p}} \ \hbox{ as }\
n\to \infty\ ,\
\hbox{ uniformly in } \ {\bf a},\ell.
\end{equation}
We shall derive (\ref{28t}) from Lemma \ref{raolemma} and Lemma \ref{aldlemma}.
As we have seen above, there exists  a separable metric $d$ on $S$,
generating the same $\sigma$-field as the Prohorov metric $\pi$,
such that (\ref{e10}) holds.
But then the limit random measure $\mu$,
which is a random variable taking values in $(S,\pi)$ ({\it i.e.},
a measurable map from the underlying probability space to $(S,{\cal B}_
{\pi})$ where ${\cal B}_{\pi}$ denotes the Borel $\sigma$-field in $S$
generated by $\pi$) can be also regarded as a random variable taking
values in $(S,d)$.  Also, $\mu$ is clearly $\sigma ({\bf X})$ measurable
and thus $(X_n,\mu) \buildrel {\cal D}\over\longrightarrow (Y_1,\mu)$
by Lemma \ref{aldlemma}.  Hence, (\ref{28t}) will follow from Lemma \ref{raolemma} (note the equivalence
of (\ref{20}) and (\ref{21}))  if we show that the class of functions
\begin{equation}\label{29t}\left\{ {{g^{{\bf a},\ell}(t,\nu)}\over {\psi ({\bf a})^p }}
\right\}
\end{equation}
defined on the product metric space $(\mathbb{R}^1\times S\ ,\ \lambda^1\times d)$
($\lambda^1$ denotes the ordinary distance on $\mathbb{R}^1$) satisfies conditions
(a),(b) of Lemma \ref{raolemma}.  Observe now that
\begin{equation}\label{30t}
\psi (a_1,\ldots,a_n)\ge \psi (a_1^*,\ldots,a_n^*)
\end{equation}
where $a_i^*$ equals either $a_i$ or 0.  (In case $(Y_n)$ is an i.i.d.\
sequence with mean 0, (\ref{30t}) follows from Jensen's inequality (see e.g.\
\cite[p.\ 153]{fe}) and the fact that, for any $H\subset \{ 1,2,\ldots,n
\}$, the conditional expectation of $\sum_{i=1}^n a_iY_i$ given $\sigma
\{Y_j\ ,\ j\in H\}$ is $\sum_{i\in H} a_iY_i$.  Since $(Y_n)$ is
a mixture of i.i.d.\ sequences with mean 0, (\ref{30t}) holds in general.)
In particular,
\begin{equation}\label{31t}
\psi (a_1,\ldots,a_n) \ge \psi(0,a_2,\ldots,a_n)
\end{equation}
and
\begin{equation}\label{32t}
\psi (a_1,\ldots,a_n) \ge \psi (a_1,0,\ldots,0) = \hbox{const } \cdot
\vert a_1\vert
\end{equation}
and thus using (\ref{e10}) we get for any $\nu \in S$, $t\in \mathbb{R}^1$ and
${\bf a} = (a_1,\ldots a_\ell)\in \mathbb{R}^\ell$,
\begin{align}\label{33t}
\big\Vert a_1t+\sum_{i=2}^\ell a_i\xi_i^{(\nu)} \big\Vert_p &
\le \vert a_1t\vert +\big\Vert \sum_{i=2}^\ell a_i\xi_i^{(\nu)}
\big\Vert_p \le \nonumber\\
&\le \vert a_1t\vert +\psi (a_2,\ldots,a_\ell)d(\nu,0) \le \hbox{ const }
\cdot \psi ({\bf a})\vert t\vert +\psi({\bf a})d(\nu,0)\nonumber\\
&\le \hbox{ const } \psi ({\bf a}) \bigl( \vert t\vert +d(\nu,0)\bigr).
\end{align}
Hence using (\ref{e10}), (\ref{31t}), (\ref{32t}) and the inequality $\vert x^p-y^p\vert \le
\vert x-y\vert \cdot p\cdot (x^{p-1}+y^{p-1})$ $(x>0,\ y>0)$ we get for any
real $t,t'$ and $\nu, \nu' \in S$
\begin{align*}&\big\vert g^{{\bf a},\ell}(t,\nu)-g^{{\bf a},\ell} (t',\nu')
\big\vert = \Big\vert \big\Vert a_1t+\sum_{i=2}^\ell
a_i\xi_i^{(\nu)} \big\Vert_p^p - \big\Vert a_1t'+ \sum_{i=2}^\ell
a_i \xi_i^{(\nu')} \big\Vert_p^p\Big\vert \cr &\le \Big\vert
\big\Vert a_1t+\sum_{i=2}^\ell a_i\xi_i^{(\nu)} \big\Vert_p
-\big\Vert a_1t'+\sum_{i=2}^\ell
a_i\xi_i^{(\nu')}\big\Vert_p\Big\vert \hbox{ const } p\psi ({\bf
a})^{p-1}\bigl( \vert t\vert + \vert t'\vert +d(\nu,0)
+d(\nu',0)\bigr)^{p-1}\cr &\le \hbox{ const}\bigl( \vert a_1\vert
\ \vert t-t'\vert +\psi (a_2,\ldots,a_n)d(\nu,\nu')\bigr)p\psi
({\bf a)}^{p-1} \bigl( \vert t\vert +\vert t'\vert
+d(\nu,0)+d(\nu',0)\bigr)^{p-1}\cr &\le \hbox{ const}\bigl( \vert
t-t'\vert +d(\nu,\nu')\bigr) p\psi({\bf a})^p \bigl( \vert t\vert
+\vert t'\vert +d(\nu,0)+d(\nu',0)\bigr)^{p-1}\cr &\le \hbox{
const}\bigl( \vert t-t'\vert +d(\nu,\nu')\bigr) p\psi ({\bf a})^p
\bigl( 2\vert t\vert +2d(\nu, 0) + \vert t-t'\vert
+d(\nu,\nu')\bigr)^ {p-1}.\cr
\end{align*}
Given $t,\nu$ and $\varep >0$, there exists a $\delta=\delta
(t,\nu,\varep)>0$ such that the last expression is $\le \varep
\psi ({\bf a})^p$ provided $\vert t-t'\vert +d(\nu,\nu') \le
\delta$ and thus the class (\ref{29t}) is locally equicontinuous
on the product metric space $(\mathbb{R}^1\times S\ ,\ \lambda^1
\times d)$.  On the other hand, (\ref{33t}) shows that the
function in (\ref{29t}) is bounded by $\hbox{const}(\vert t\vert
+d(\nu,0))^p\le \hbox{const }2^p (\vert t\vert^p +d(\nu,0)^p)$.
Now, using $(X_n,\mu)\buildrel {\cal D} \over\longrightarrow
(Y_1,\mu)$, the uniform integrability of $\vert X_n \vert^p$ and
$Ed(\mu,0)^p <+\infty$ (see (\ref{ed})) we get
$$
E\bigl( \vert X_n\vert^p +d(\mu,0)^p\bigr) \longrightarrow
E\bigl( \vert Y_1\vert^p +d(\mu, 0)^p\bigr).
$$
Thus the class
(\ref{29t}) satisfies also condition (b) of Lemma \ref{raolemma}.
We thus proved relation (\ref{28t}) and thus also (\ref{25t})
whence it follows (note again that the right side of (\ref{25t})
equals 1) that
\begin{align}\label{34t}
&\psi ({\bf a})^{-1}\Vert a_1X_n+a_2Y_2+\cdots +a_\ell Y_\ell \Vert_p
 \nonumber \\
& \phantom{999999999999999} \longrightarrow \psi ({\bf a})^{-1}\Vert a_1Y_1+a_2Y_2+\cdots +a_\ell Y_\ell\Vert_p
\ \hbox{ as }\ n\to \infty \end{align}
unformly in $\ell, {\bf a}$.  Hence we can choose $n_1$ so large that
$$
\big\vert \ \Vert a_1X_{n_1} +a_2Y_2 +\cdots +a_\ell Y_\ell \Vert_p
-\Vert a_1Y_1+a_2Y_2+\cdots +a_\ell Y_\ell \Vert_p \big\vert
\le {\varep \over 2} \psi (a_1,\ldots,a_\ell)$$
for every $\ell,{\bf a}$.
This completes the first induction step.

Assume now that $n_1,\ldots ,n_{k-1}$ have already been chosen.
Exactly in the same way as we proved (\ref{34t}), it follows that for $\ell >k$
\begin{align*}&\psi ({\bf a})^{-1}\Vert a_1X_{n_1} +\cdots + a_{k-1}X_{n_{k-1}}
+a_k X_n +a_{k+1}Y_{k+1}+\cdots + a_\ell Y_\ell \Vert_p \cr
&\longrightarrow \psi ({\bf a})^{-1}\Vert a_1X_{n_1}+\cdots + a_{k-1}
X_{n_{k-1}} +a_kY_k+\cdots +a_\ell Y_\ell \Vert_p \ \hbox{ as }\
n\to \infty\cr
\end{align*}
uniformly in {\bf a} and $\ell$.  Hence we can choose $n_k$ so large
that $n_k>n_{k-1}$ and
\begin{align*}&\big\vert \ \Vert a_1X_{n_1}+\cdots +a_{k-1}X_{n_{k-1}}
+a_kX_{n_k} +a_{k+1}Y_{k+1} +\cdots + a_\ell Y_\ell\Vert_p\cr
&\qquad - \Vert a_1X_{n_1}+\cdots + a_{k-1}X_{n_{k-1}} +a_kY_k +\cdots
+ a_\ell Y_\ell \Vert_p\big\vert \le {{\varep}\over{2^k}} \psi
(a_1,\ldots ,a_\ell) \cr
\end{align*}
for every $(a_1,\ldots,a_\ell)\in \mathbb{R}^\ell$ and $\ell >k$.  This completes
the $k$-th induction step; the so constructed sequence $(n_k)$
obviously satisfies
$$\big\vert \ \Vert a_1X_{n_1}+\cdots +a_\ell X_{n_\ell} \Vert_p
- \Vert a_1Y_1 +\cdots + a_\ell Y_\ell \Vert_p\vert \le \varep
\psi (a_1,\ldots ,a_\ell)$$ for every $\ell \ge 1$  and
$(a_1,\ldots, a_\ell)\in \mathbb{R}^\ell$. The last relation is
equivalent to (\ref{24t}) and thus the sufficiency of
(\ref{maincond}) in Theorem \ref{th2} is proved.

\medskip
We now turn to the proof of necessity of (\ref{maincond}) in Theorem \ref{th2}.
Assume that $(X_n)$ is equivalent to the
unit vector basis of $\ell^2$; then for any increasing sequence $(m_k)$ of integers
we have
$$
\left\|\frac{1}{\sqrt{N}}  \sum_{k=1}^N X_{m_k}\right\|_p=O(1)
$$
and thus by the Markov inequality we have for any $A\subset \Omega$ with $P(A)>0$,
\begin{equation} \label{bp}
P_A \left\{ \left|\frac{1}{\sqrt{N}} \sum_{k=1}^N X_{m_k}\right| \ge T \right\}\le 1/2
\qquad \text{for} \ T\ge T_0,\, N\ge 1
\end{equation}
where $T_0$ depends on $A$ and the sequence $(X_n)$. We show first that
\begin{equation}\label{l2int}
 \int_{-\infty}^\infty x^2 d\mu (x)<\infty \qquad \text{a.s.}
\end{equation}
Let $F_\bullet (x)$ denote the random distribution function corresponding to
$\mu$ and assume indirectly that there exists a set $B \subset \Omega$ with $P(B) > 0$
such that
\begin{equation}\label{3.2}
\lim_{t \to \infty} \int_{|x| < t} x^2 dF_\bullet(x) = +\infty
\quad \text{on } B.
\end{equation}
By Egorov's theorem there exists a set $B^* \subset B$ with $P(B^*) \geq
P(B)/2$ such that on $B^*$ (\ref{3.2})
holds uniformly, i.e.\ there exists  a positive, nondecreasing, nonrandom function $K_t \to +\infty$ such that
\begin{equation}\label {3.9a}
\int_{|x| < t} x^2 dF_\bullet (x) \geq K_t  \ \text{ on } \ B^* .
\end{equation}
Also,
\begin{equation}\label {3.9b}
 \int_{|x|\ge t} dF_\bullet(x) \longrightarrow 0 \quad  \text{a.s.\  as}  \ t \to
\infty
\end{equation}
and thus we can choose a set $B^{**} \subset B^*$ with
$P(B^{**}) \geq P(B^*) / 2$ such that on $B^{**}$ relation (\ref{3.9b})
holds uniformly, i.e.\ there exists a positive, nonincreasing, nonrandom
function $\tilde{\ve}_t \to 0$ such that
\begin{equation}\label{3.3}
\int_{|x|\ge t} dF_\bullet(x) \leq \tilde{\ve}_t \ \text{ on } \
B^{**}.
\end{equation}
We show that there exists a subsequence $(X_{m_k})$ of
$(X_n)$ such that (\ref{bp}) fails for $A= B^{**}$. Since our
argument will involve the sequence $(X_n)$ only on the set
$B^{**}$ and on $B^{**}$ $(X_n)$ satisfies the assumptions of Theorem \ref{th2}
with the same $\mu$ and with $\|X_n\|_p=1$ replaced by $\|X_n\|_p=O(1)$ (which is all
we need for the rest of the proof), in the sequel we can assume, without loss of
generality, that $B^{**} = \Omega$. That is, we may assume that (\ref{3.9a}),
(\ref{3.3}) hold on the whole probability space.

Let $C$ be an arbitrary set in the probability space with $P(C)
> 0$. Integrating (\ref{3.9a}), (\ref{3.3}) on $C$ and using  (\ref{(4)}) and Lemma \ref{lemma4} we get
\begin{equation}\label{3.5}
\int_{|x| < t} x^2 dF_C(x) \geq K_t, \qquad
 \int_{|x|\ge t} dF_C (x) \leq \tilde{\ve}_t \quad \text{for} \ t, -t \in H
\end{equation}
where $H$ denotes the set of continuity points of $F_C$.
Since the integrals in (\ref{3.5}) are monotone functions of
$t$ and ${\mathbb R}\setminus H$ is countable,  (\ref{3.5}) remains valid with $K_{t/2}$ resp.
$\tilde{\ve}_{t/2}$ if we drop the assumption $t, -t \in H$. Thus, keeping the original notation, in the sequel
we can assume that (\ref{3.5}) holds for all $t>0$. Choose now $t_0$ so
large that $\tilde{\ve}_{t_0} \leq 1/16$ and then
choose $t_1>t_0$ so large that
$$
K^{1/2}_t \ge 4t_0 \ \text{ for } \ t \geq t_1.
$$
Then for $t \geq t_1$ we have, using
(\ref{3.5}),
\begin{align*}
&\biggl| \int_{|x| < t} x dF_C(x) \biggr| \leq t_0 + \int\limits_{t_0
\leq |x| < t} |x| dF_C(x) \\
&\leq t_0 + \biggl( \int\limits_{|x| \geq t_0} dF_C(x) \biggr)^{1/2}
\biggl( \int\limits_{|x| < t} x^2 dF_C(x) \biggr)^{1/2} \\
&\leq t_0 + {1 \over 4 } \biggl( \int\limits_{|x| < t} x^2 dF_C(x)
\biggr)^{1/2} \leq {1\over 2} \biggl(\int\limits_{|x| < t} x^2
dF_C(x) \biggr)^{1/2}
\end{align*}
and thus for any $C \subset \Omega$ with $P(C) >0$ we have
\begin{equation}\label{3.6}
\int_{|x|< t} x^2 dF_C(x) - 2 \biggl( \int_{|x|< t}
xdF_C(x)\biggr)^2 \geq {1\over 2} K_t, \
\quad t \geq t_1.
\end{equation}
Let now $(\ve_n)$ tend to $0$ so rapidly that
\begin{equation} \label{doublestar}
\sum_{j=a_k+1}^\infty \ve_j \leq 2^{-k}.
\end{equation}
Let $a_k=[\log k +1]$ $(k=1, 2, \ldots)$.
By Lemma \ref{lemma1} there exists a
subsequence $(X_{m_k})$ and a sequence $(Y_k)$ of discrete
r.v.'s such that (\ref{1}) holds and for each $k \ge 1$ the atoms of
the finite $\sigma$-field $\sigma \{ Y_1, \dots, Y_{a_k} \}$ can
be divided into two classes $\Gamma_1$ and $\Gamma_2$ such that
\begin{equation}\label{3.9}
 \sum_{B \in \Gamma_1} P(B) \leq \ve_{a_k+1} \leq 2^{-k}
\end{equation}
and for each $B \in \Gamma_2$ there exist $P_B$-independent r.v.'s
$Z^{(B)}_{a_k + 1}, \dots, Z^{(B)}_k$ defined on $B$ with common distribution
$F_B$ such that
\begin{equation}\label{3.10}
P_B \bigl( |Y_j - Z^{(B)}_j | \geq 2^{-k} \bigr) \leq 2^{-k}
\quad (j = a_k + 1, \dots, k).
\end{equation}
Set
\begin{align*}
S^{(B)}_{a_k, k} &= \sum^k_{j = a_k + 1} Z^{(B)}_j, \qquad B \in
\Gamma_2 \\
\overline S_{a_k, k} &= \sum_{B \in \Gamma_2} S^{(B)}_{a_k, k} I(B),
\end{align*}
where $I(B)$ denotes the indicator function of $B$.
By (\ref{3.10}) and $k 2^{-k}\le 1$,
$$
P_B \biggl( \biggl| \sum^k_{j = a_k + 1} Y_j - \sum^k_{j = a_k +
1} Z^{(B)}_j \biggr| \geq 1 \biggr) \leq k2^{-k}, \qquad B \in
\Gamma_2
$$
and thus using (\ref{3.9}) we get
\begin{equation}\label{3.12}
P \biggl( \biggl| \sum^k_{j = a_k + 1} Y_j - \overline S_{a_k,
k} \biggr| \geq 1 \biggr) \leq (k+1) 2^{-k}.
\end{equation}
Since $\|X_n\|_1=O(1)$, by the triangular inequality and the Markov inequality we have
\begin{align*}
&P \biggl( \biggl| \sum^{a_k}_{j = 1} X_{m_j} \biggr| \geq a_k
k^{1/4} \biggr) \leq \sum_{j=1}^{a_k} P \bigl(
|X_{m_j}| \geq k^{1/4} \bigr)\\
& \leq \text{const} \, (\log k+1) k^{-1/4} =: \delta_k
\end{align*}
which, together with (\ref{3.12}), (\ref{1}) and (\ref{doublestar}), yields
\begin{align}\label{3.13}
&P \biggl( \biggl| \sum^k_{j = 1} X_{m_j} - \overline S_{a_k, k}
\biggr| \geq 3a_k k^{1/4} \biggr)
\le P \biggl( \biggl| \sum^{a_k}_{j = 1} X_{m_j} \biggr| \geq a_k
k^{1/4} \biggr) \nonumber \\
& + P \biggl( \sum^k_{j = a_k+ 1} |X_{m_j} -Y_j| \geq 1 \biggr)
+ P \biggl( \biggl| \sum^k_{j = a_k+ 1} Y_j - \overline S_{a_k, k}
\biggr| \geq 1 \biggr)  \nonumber \\
&\leq \delta_k +(k+2) 2^{-k}.
\end{align}
Applying Lemma \ref{lemma2} to the i.i.d.\ sequence $\{ Z^{(B)}_j, a_k + 1
\leq j \leq k \}$ and using (\ref{3.6}) with $C=B$, $a_k\le k/2$
and the monotonicity of $K_t$ we get for any $T\ge 2$,
\begin{align*}
& P_B \left(\biggl| {S^{(B)}_{a_k, k} \over \sqrt k} \biggr|
\leq T \right)
\leq P_B \left( \left|{S^{(B)}_{a_k, k} \over \sqrt{k -
a_k}}\right| \leq 2T \right) \leq
\text{\rm const} \cdot 2T K^{-1/2}_{2T \sqrt {k-a_k}} \\
&\le
\text{\rm const} \cdot T K^{-1/2}_{2T\sqrt {k/2}}
\end{align*}
where the constants are absolute. Thus using (\ref{3.9}) it follows that
\begin{equation}\label{3.14}
P \biggl( \biggl| {\overline S_{a_k, k} \over \sqrt k} \biggr|
\leq T  \biggr) \leq \text{\rm const} \cdot TK^{-1/2}_{\sqrt {k/2}} + 2^{-k}.
\end{equation}
Using (\ref{3.13}), (\ref{3.14})  and $a_k \leq \log k + 1$  it follows that
\begin{align*}
&P\left( \left| \frac{1}{\sqrt{k}}\sum_{j=1}^k X_{m_j}\right| \le T\right) \le
P \left( \left| {\overline S_{a_k, k} \over \sqrt k} \right|\leq T +3a_k k^{-1/4} \right) + (k+2) 2^{-k} + \delta_k\\
& \le \text{\rm const} \, \left(T+3a_kk^{-1/4}\right)  K^{-1/2}_{\sqrt {k/2}} + (k+2) 2^{-k} + \delta_k
\longrightarrow 0\qquad \text{as} \ k\to\infty
\end{align*}
for any fixed $T\ge 2$ which clearly contradicts to (\ref{bp}) with $A=\Omega$. This completes the proof of (\ref{l2int}).

Since $X_n\longrightarrow 0$ weakly in $L^p$, we have $\int_{-\infty}^\infty x d\mu(x)=0$ a.s., and
thus the random measure $\mu$ has mean zero and finite variance with probability one. Thus the subsequence principle, specialized to
the central limit theorem (see e.g.\ \cite[Theorem 3]{bepe} and the remark following it)
there exists a subsequence $(X_{n_k})$ such that
\begin{equation}\label{G}
N^{-1/2} \sum_{k=1}^N X_{n_k} \overset{\cal D}{\longrightarrow} Y\zeta
\end{equation}
where
$Y=\left(\int_{-\infty}^\infty x^2 d\mu (x)\right)^{1/2}$, $\zeta$
is a standard normal variable and $Y$ and $\zeta$ are independent.  Note that (\ref{G}) holds in distribution, but by a well known result of Skorokhod
(see e.g.\ \cite{bill}, p.\ 70) there exist r.v.'s $W_N$, $W$ $(N=1, 2, \ldots)$ such that  $W_N$ has the same
distribution as $N^{-1/2} \sum_{k=1}^N X_{n_k}$ in (\ref{G}), $W$ has the same distribution as $Y\zeta$ and $W_N\longrightarrow W$
a.s. Thus (\ref{G}) and Fatou's lemma imply
\begin{equation}\label{ut}
\|Y \zeta\|_p \le \liminf_{N\to\infty} \| N^{-1/2} \sum_{k=1}^N X_{n_k}\|_p<\infty
\end{equation}
where the second inequality follows from the equivalence of $(X_n)$ to the unit vector basis of
$\ell^2$,
assumed at the beginning of the proof. Since $Y$ and $\zeta$ are independent, (\ref{ut}) implies $E|Y|^p<\infty$,
i.e.\ (\ref{maincond}) holds, completing the proof of the converse part of Theorem \ref{th2}.

\bigskip\bigskip\medskip\noindent
{\bf Acknowledgement.} The authors are indebted to the referee for his/her valuable remarks and suggestions leading to a considerable improvement of the presentation.

\end{document}